\documentclass[11pt]{article}
 \usepackage{amsmath, amsthm, amssymb, bbm, setspace,bigints}
  \usepackage{algorithm2e}
\RestyleAlgo{boxruled}
\LinesNumbered
 \usepackage[margin=1 in]{geometry}
\usepackage{graphicx}				
\usepackage{amssymb,color}					
\usepackage{amsmath}
\usepackage{epstopdf}
\usepackage{bbm}
\usepackage{mathtools}
\usepackage{authblk}
\usepackage{mathrsfs}
\newcommand{\argmax}{\mathop{\rm argmax~}}

\usepackage{graphicx,array,booktabs,bbm,multicol,colortbl,tabularx}
\usepackage{caption,subcaption}

\newcommand{\argmin}{\mathop{\rm argmin~}}

\usepackage{mathtools}

\DeclarePairedDelimiterX{\infdivx}[2]{(}{)}{%
  #1\;\delimsize\|\;#2%
}

\newcommand{\qhat} {\widehat{q}}
\newcommand{\qtil} {\widetilde{q}}

\newcommand{\norm}[1]{\left\Vert#1\right\Vert}

\def\m{\mathcal}
\def\mb{\mathbb}

\newcommand \bbP{\mathbb{P}}
\newcommand \bbE{\mathbb{E}}

\newtheorem{theorem}{Theorem}[section]
\newtheorem{lemma}[theorem]{Lemma}

\newtheorem{corollary}[theorem]{Corollary}

\begin{document}

\title{On Statistical Optimality of Variational Bayes\footnote{To appear in AISTATS 2018}}

\author[1]{Debdeep Pati\thanks{debdeep@stat.tamu.edu}}
\author[1]{Anirban Bhattacharya\thanks{anirbanb@stat.tamu.edu}}
\author[2]{Yun Yang \thanks{yyang@stat.fsu.edu}}
\affil[2]{Department of Statistics, Florida State University}
\affil[1]{Department of Statistics, Texas A\&M University}
\date{\vspace{-2em}}
\maketitle 

\begin{abstract}
The article  addresses a long-standing open problem on the justification of using variational Bayes methods for parameter estimation.  We provide general conditions for obtaining optimal risk bounds for point estimates acquired from mean-field variational Bayesian inference. The conditions pertain to the existence of certain test functions for the distance metric on the parameter space and minimal assumptions on the prior.  A general recipe for verification of the conditions is outlined which is  broadly applicable to existing Bayesian models with or without latent variables. As illustrations, 
specific applications to Latent Dirichlet Allocation and Gaussian mixture models are discussed. 
\end{abstract}

\section{Introduction}

Variational inference \cite{jordan1999introduction,bishop2006pattern,wainwright2008graphical} is now a well-established tool to approximate intractable posterior distributions in hierarchical multi-layered Bayesian models. The traditional Markov chain Monte Carlo (MCMC; \cite{gelfand1990sampling}) approach of approximating distributions with intractable normalizing constants draws (correlated) samples according to a discrete-time Markov chain whose stationary distribution is the target distribution. Despite their success and popularity, MCMC methods can be slow to converge and lack scalability in big data problems and/or problems involving very many latent variables, which has fueled search for alternatives. 

In contrast to the sampling approach of MCMC, variational inference approaches the problem from an optimization viewpoint. First, a class of analytically tractable distributions, referred to as the {\em variational family}, is identified for the problem at hand. For example, in {\em mean-field} approximation, the set of parameters and latent variables is divided into blocks and the variational distribution is assumed to be independent across blocks. The distribution in the variational family closest to the target distribution relative to the Kullback--Leibler (KL) divergence is then used as a proxy to the target. Implementation-wise, the above optimization is commonly solved using coordinate descent or alternating minimization. A comprehensive review of various aspects of variational inference can be found in the recent article by \cite{blei2017variational}.

Variational inference has arguably found its most potent applications in latent variable models such as mixture models, hidden Markov models, graphical models, topic models, and neural networks; see \cite{mackay1997ensemble,attias2000variational,wainwright2008graphical,humphreys2000approximate,corduneanu2001variational,ueda2002bayesian,blei2012probabilistic,blei2003latent,graves2011practical} for a flavor of this enormous literature. Due to the fast convergence properties of the variational objective, variational inference algorithms are typically orders of magnitude faster in big data problems compared to MCMC approaches \cite{kushner1997stochastic,kingma2014adam,ahmed2012scalable}. However, in spite of their tremendous empirical success, a general statistical theory qualifying the statistical properties of a variational solution is lacking. Existing results operate in a case-by-case manner, often directly analyzing the specific iterative algorithm to characterize properties of its limit; see Section 5.2 of \cite{blei2017variational} for a comprehensive review. These analyses typically require sufficient tractability of the successive steps of the iterative algorithm, and can be difficult to adapt to minor changes in the prior. More recently \cite{alquier2016properties,yang2017variational,alquier2017concentration} modified the objective function to introduce an inverse-temperature parameter, and obtained general guarantees for the variational solution under this modified objective function.

In this article, we aim to address the general question as to whether point estimates obtained from usual variational approximations share the same statistical accuracy as those from the actual posterior. We clarify at the onset that we operate in a frequentist setting assuming the existence of a true data generating parameter. Our novel contribution is to relate the Bayes risk relative to the variational solution for a general distance metric (defined on the parameter space) to (i) the existence of certain test functions for testing the true parameter against complements of its neighborhood under the error metric, and (ii) the size of the variational objective function. As an important consequence of such Bayes risk bounds, risk bounds for variational point estimates can be readily derived when the distance metric is convex. If the risk of the variational point estimate coincides with the contraction rate of the exact posterior at the true parameter, it can be argued that there is no loss of statistical efficiency, at least asymptotically, in using the variational approximation. 

We identify a number of popularly used models where the conditions can be satisfied and the variational point estimate attains the minimax rate. Since variational Bayes is primarily used for point estimation, our theory suggests that variational Bayes successfully achieves its desiderata. As vignettes, we present two non-trivial examples in the form of density estimation in Latent Dirichlet Allocation (LDA) for topic modeling, and estimating component specific parameters in Gaussian mixture models.

\section{Background}

In this section, we introduce notations and offer some background to set up our theoretical results. Let $h(p\,||\,q) =[\int(\sqrt{p} - \sqrt{q})^2d\mu]^{1/2}$ and $D(p\,||\,q) = \int p\log (p/q)d\mu$ denote the Hellinger distance and Kullback--Leibler divergence, respectively, between two probability density functions $p$ and $q$ relative to a a common dominating measure $\mu$. Also denote $V(p\,||\,q) = \int p \log^2(p/q) d\mu$. For a set $A$, we use $I_A$ to denote its indicator function. For any vector $\mu$ and positive semidefinite matrix $\Sigma$, we use $\m N(\mu, \Sigma)$ to denote the normal distribution with mean $\mu$ and covariance matrix $\Sigma$, and use $\m N(\theta;\, \mu, \Sigma)$ to denote its pdf at $\theta$. We use w.h.p. to abbreviate ``with high probability", when the probability is evident from the context. Throughout, $C$ denotes a constant independent of everything else whose value may change from one line to the other. We write $a \lesssim b$ to denote $a \le C b$ for some constant $C> 0$. Similarly, $a \gtrsim b$. 

Suppose we have $n$ observations $Y_1, \ldots, Y_n$ and a probabilistic model $\mb P_{\theta}^{(n)}$ for the joint distribution of the data $Y^n=(Y_1,\ldots,Y_n)$, with a density $p_{\theta}^{(n)}$ relative to the Lebesgue measure. Here, $\theta$ is the unknown parameter to be estimated from the data which lives in some parameter space $\Theta$. Our formulation does not require the $Y_is$ to be identically distributed or even independent. In a Bayesian setup, uncertainty regarding the parameter is quantified through a prior distribution $p_{\theta}$ on $\Theta$, which upon observing the data is updated to the posterior distribution using Bayes' theorem: 
\begin{align}\label{eq:post}
p(\theta\in B\,|\,Y^n) = \frac{\int_{B} \big[p(Y^n\,|\,\theta)\big]\, p_\theta(\theta)\, d\theta}{\int_{\Theta} \big[p(Y^n\,|\,\theta)\big]\, p_\theta(\theta)\, d\theta},
\end{align}
for any measurable subset $B \subset \Theta$ and $p(Y^n\,|\,\theta) := p_{\theta}(Y^n)$.  

In a wide variety of practical problems, the likelihood function $p(Y^n \,|\, \theta)$ may be intractable or difficult to analyze directly. For example, in a 2-component Gaussian mixture model, $p(Y^n \,|\, \theta)$ is a combinatorial sum over $2^n$ terms. A widely used trick in such situations is to introduce latent variables to simplify the conditional likelihood. Specifically, assume one can decompose
\begin{align}\label{eq:latvar}
p(Y^n \,|\, \theta) = \sum_{s^n} p(Y^n \,|\, S^n = s^n, \mu) \, \pi_{s^n}, 
\end{align}
where $S^n = (S_1, \ldots, S_n)$ denotes a collection of discrete latent variables, with $S_i\in\{1,2,\ldots,K\}$ the latent variable for the $i$th observation. We have assumed the parameter $\theta$ can be decomposed as $\theta = (\mu, \pi)$, with $p(Y^n \,|\, S^n = s^n, \theta) = p(Y^n\,|\,S^n = s^n, \mu)$ and 
and $\pi_{s^n} :\,= p(S^n = s^n \,|\, \theta)$ denotes the probability of the latent vector taking on the value $s^n$. In the 2-component mixture model example, $S_i \in \{1, 2\}$ denotes the latent membership indicator for the $i$th observation.  We assume discrete latent variables for notational convenience and note that our results generalize to continuous latent variables in a straightforward fashion. 

Let $Z^n = (\theta, S^n)$. The augmented posterior $p(Z^n \,|\, Y^n)$ assumes the form 
\begin{align}\label{eq:aug_post}
p(\theta, s^n \,|\, Y^n) \propto p(Y^n \,|\, \mu, s^n) \, \pi_{s^n} \, p_{\theta}(\theta), 
\end{align}
whose constituent terms are typically more tractable. 
Variational inference in this setup proceeds by first identifying a {\em variational family} comprising of distributions $\Gamma$ on $Z^n$ and finding the closest member in this family to $p(Z^n\,|\,Y^n)$ relative to the KL divergence 
\begin{align}\label{eq:VB_obj1}
\qhat_{Z^n}&:\, = \argmin_{q_{Z^n} \in \Gamma} D\big[\,q_{Z^n}(\cdot)\,\big|\big|\, p(\cdot\,|\, Y^n)\,\big] = \argmax_{q_{Z^n} \in \Gamma} \m L(q_{Z^n}),
\end{align}
where 
\begin{align}\label{eq:ELBO}
\m L(q_{Z^n}) = \int q_{Z^n}(z^n)\, \log \frac{p(Y^n\,|\,z^n)\, p_{Z^n}(z^n)}{q_{Z^n}(z^n)}\, dz^n
\end{align}
is the {\em evidence lower bound} (ELBO) which gives a lower bound to the log marginal likelihood $\log p(Y^n)$. If $\Gamma$ is completely unrestricted, then $\qhat_{Z^n}$ coincides with the posterior distribution $p(Z^n\,|\,Y^n)$. In practice, the choice of the variational family $\Gamma$ is dictated by a trade-off between flexibility and computational tractability. For example, in {\em mean-field} variational approximation, it is common to assume independence among the parameters and the latent variables in the variational family, whence $q_{Z^n}$ decomposes as 
\begin{align}\label{eq:MF}
q_{Z^n}(z^n) = q_{\theta}(\theta) \otimes q_{S^n}(s^n). 
\end{align}
{\em We shall assume the mean-field decomposition \eqref{eq:MF} throughout}, so that the minimizer $\qhat_{Z^n}$ in \eqref{eq:VB_obj1} necessarily is of the form 
$\qhat_{Z^n} = \qhat_{\theta} \otimes \qhat_{S^n}.
$
In many situations, the constituent terms may be further decomposed as $q_{S^n}(s^n) = \otimes_{i=1}^n q_{S_i}(s_i)$ and $q_{\theta}(\theta) = q_{\mu}(\mu) \otimes q_{\pi}(\pi)$. Such a decomposition may be either due to computational reasons or implied by the conditional independence structure of the model.

Under the mean-field decomposition \eqref{eq:MF} and using \eqref{eq:aug_post}, we have, after some simplification, 
\begin{align}\label{eq:KL_decomp}
& D\big[\,q_{Z^n}(\cdot)\,\big|\big|\, p(\cdot\,|\, Y^n)\,\big] = D(q_{\theta}\,\big|\big|\,p_{\theta}) \notag \\
& - \int_{\Theta} \underbrace{ \bigg[ \sum_{s^n} q_{S^n}(s^n) \, \log \frac{p(Y^n\,|\,\mu, s^n) \, \pi_{s^n}}{q_{S^n}(s^n)} \bigg] }_{\widetilde{\ell}_n(\theta)} \, q_{\theta}(d\theta).
\end{align}
The quantity $\widetilde{\ell}_n(\theta)$ is an approximation to the log likelihood $\ell_n(\theta) :\,= \log p(Y^n \mid \theta)$ in terms of the latent variables. To see this, multiply and divide the right hand side of \eqref{eq:latvar} by $q_{S^n}(s^n)$ and apply Jensen's inequality to the concave function $x \mapsto \log x$ to conclude that $\ell_n(\theta) \ge \widetilde{\ell}_n(\theta)$. Replacing $\widetilde{\ell}_n(\theta)$ by $\ell_n(\theta)$ in \eqref{eq:KL_decomp1} and adjusting the error term, we obtain
\begin{align}\label{eq:KL_decomp1}
& D\big[\,q_{Z^n}(\cdot)\,\big|\big|\, p(\cdot\,|\, Y^n)\,\big] = D(q_{\theta}\,\big|\big|\,p_{\theta}) - \int_{\Theta} \ell_n(\theta) q_{\theta}(d\theta) \notag \\
& + \underbrace{ \int_{\Theta}  \big[ \ell_n(\theta) - \widetilde{\ell}_n(\theta) \big]  \, q_{\theta}(d\theta)}_{\Delta(q_{\theta}, q_{S^n})}.
\end{align}
The quantity $\Delta$ is an average error due to the likelihood approximation and is clearly nonnegative. In the specific situation where no latent variables are present in the model, $\Delta \equiv 0$. However, in general, $\Delta$ is a strictly positive quantity. 

\section{Variational risk bounds}

We are now prepared to state a general Bayes risk bound for the variational distribution. We shall operate in a frequentist framework and assume the existence of a true data generating parameter $\theta^\ast \in \Theta$. In other words, we assume that the data $Y^n$ is distributed according to $\mb P_{\theta^\ast}^{(n)}$. Let $d(\cdot, \cdot)$ be a distance metric on the parameter space $\Theta$ which quantifies distance between two putative parameter values. For example, if $\Theta$ designates a space of densities so that each $\theta$ can be identified with a density function, $d$ can be chosen as the Hellinger or total variation distance. We are interested in obtaining bounds on 
$$
\int d^{\kappa}(\theta, \theta^\ast) \,\qhat_{\theta}(\theta) 
$$
for some $\kappa > 0$, that hold with high probability with respect to $\mb P_{\theta^\ast}^{(n)}$. In particular, if $d^{\kappa}$ is convex in its first argument, then by Jensen's inequality, 
$$d^{\kappa}(\widehat{\theta}_{VB}, \theta^\ast) \le \int d^{\kappa}(\theta, \theta^\ast) \, \qhat_{\theta}(\theta),$$ where $\widehat{\theta}_{VB} = \int \theta \,\qhat_{\theta}(d\theta)$ is the mean of the variational distribution and a surrogate for the posterior mean. We are specifically interested in obtaining sufficient conditions for the variational point estimate $\widehat{\theta}_{VB}$ to contract at the same rate as posterior mean. Since variational approaches are overwhelmingly used for rapidly obtaining point estimates, such a result will indicate that at least in terms of rates of convergence, there is no loss of statistical accuracy in using a variational approximation. Moreover, the negative result from \cite{wang2005inadequacy} shows that the spread of the variational distribution $\qhat_{\theta}$ is typically ``too small" compared with that of the sampling distribution of the maximum likelihood estimator. This fact combined with the Bernstein von-Mises theorem (Chapter 10 of \cite{van2000asymptotic}) implies the inadequacy of using $\qhat_{\theta}$ for approximating the true posterior distribution, and a rate optimal variational point estimator is the best one can hope for in general.

Define $\ell_n(\theta, \theta^\ast) :\, = \ell_n(\theta) - \ell_n(\theta^\ast)$ to be the log-likelihood ratio between $\theta$ and $\theta^\ast$. We can replace $\ell_n(\theta)$ with $\ell_n(\theta, \theta^\ast)$ inside the integrand in \eqref{eq:KL_decomp1} without affecting the minimization problem - this is done for purely theoretical reasons to harness the structure of the log-likelihood ratio. Let us call the equivalent objective function $\Omega$, so that 
\begin{align}\label{eq:V}
& \Omega(q_{\theta}, q_{S^n}) : =\,  \notag \\
& - \int_{\Theta} \ell_n(\theta, \theta^\ast) \, q_{\theta}(d\theta) + \Delta(q_{\theta}, q_{S^n}) + D(q_{\theta}\,\big|\big|\,p_{\theta}).
\end{align}

We are now ready to state the main assumption. \\
\noindent {\bf Assumption T:  (existence of tests)} Let $\varepsilon_n$ be a sequence satisfying $\varepsilon_n \to 0$ and $n \varepsilon_n^{\kappa} \to \infty$ for some $\kappa \ge 2$. Let $\phi_n \equiv \phi_n(Y^n) \in (0, 1)$ be a sequence of test functions for testing 
$$
H_0: \theta = \theta^\ast \, \text{versus.} \, H_1: d(\theta, \theta^\ast)  > \varepsilon_n
$$
with type-I and II error rates satisfying
\begin{align*}
& \mb E_{\theta^\ast} [\phi_n] \le e^{-2 n \varepsilon_n^{\kappa}}, \\
& \mb E_{\theta} [1 - \phi_n] \le e^{-C n d^{\kappa}(\theta, \theta^\ast)},
\end{align*}
for any $\theta \in \Theta$ with $d(\theta, \theta^\ast) > \varepsilon_n$, where $\mb E_{\theta}$ denotes an expectation with respect to $\mb P_{\theta}^{(n)}$. While $\kappa = 2$ appears naturally in most problems, we provide a non-standard example (estimating component specific means in a mixture model) in Section 4 with $\kappa = 4$. 

In models satisfying the monotone likelihood ratio property \cite{casella2002statistical} such
as exponential families, one can construct such tests (with $\kappa = 2$) from the generalized likelihood ratio test (GLRT) statistic when $d$ corresponds to the Euclidean metric on the natural parameter. 
A general recipe \cite{ghosal2000} to construct such tests when $\Theta$ is compact relative to $d$ is to (i) construct an $\varepsilon_n/2$-net $\m N = \{\theta_1, \ldots,  \theta_{N}\}$ such that for any $\theta$ with $d(\theta, \theta^\ast) > \varepsilon_n$, there exists $\theta_j \in \m N$ with $d(\theta, \theta_j) < \varepsilon_n/2$, (ii) construct a test $\phi_{n,j}$ for $H_0: \theta = \theta^\ast$ versus $H_1: \theta = \theta_j$ with type-I and II error rates as in Assumption {\bf T}, and (iii) set $\phi_n = \max_{1\le j \le N} \phi_{n, j}$. 
The type-II error of $\phi_n$ retains the same upper bound, while the type-I error can be bounded by $N \, e^{-2 n \varepsilon_n^{\kappa}}$. Since $N$ can be further bounded by $N(\Theta, \varepsilon_n/2, d)$, the covering number of $\Theta$ by $d$-balls of radius $\varepsilon_n/2$, it suffices to show that $N(\Theta, \varepsilon_n/2, d) \lesssim e^{n \varepsilon_n^{\kappa}}$. For example, when $\Theta$ is a compact subset of $\mathbb{R}^d$, then $N(\Theta, \varepsilon_n/2, d) \lesssim \varepsilon_n^{-d} \lesssim e^{n \varepsilon_n^{2}}$ as long as $\varepsilon_n \lesssim \sqrt{\log n/n}$. More generally, if $\Theta$ is a space of densities and $d$ the Hellinger/$L_1$ metric, then construction of the point-by-point tests in (i) (with $\kappa = 2$) from the LRT statistics follows from the classical Birg\'{e}-Lecam testing theory \cite{birge1983approximation,lecam1973convergence}; see also \cite{ghosal2000}. 

We are now ready to state our first theorem on the variational risk bound, which relates the Bayes risk under the variational solution to the size of the objective function $\Omega$. 
\begin{theorem}\label{theorem:main}
Suppose there exists a sequence of test functions $\{\phi_n\}$ for the metric $d$ satisfying Assumption {\bf T}. Then, it holds with $\mb P_{\theta^\ast}^{(n)}$ probability at least $(1 - e^{- C n \varepsilon_n^{\kappa}})$ that 
$$
\int_{\Theta} d^{\kappa}(\theta, \theta^\ast) \, \qhat_{\theta}(d\theta) \le \frac{\Omega(q_{\theta}, q_{S^n})}{n} + \varepsilon_n^{\kappa},
$$
for any $q_{\theta} \ll p_{\theta}$ and any probability distribution $q_{S^n}$ on $S^n$ which is nowhere zero. 
\end{theorem}
Theorem \ref{theorem:main} implies that minimizing the Bayes risk is equivalent to minimizing the objective function $\Omega$ in \eqref{eq:V}. \cite{yang2017variational} obtained a similar result for a modified variational objective function and $d$ being limited to R{\'e}nyi divergence measures. Theorem \ref{theorem:main} instead allows any metric $d$ as long as the testing condition in Assumption {\bf T} can be satisfied. 

A detailed proof of Theorem \ref{theorem:main} (as well as proofs of other results) is provided in the Section \ref{appendix}. We sketch some of the key steps to highlight the main features of our argument. Our first key step is to show using the testing assumption {\bf T} that w.h.p. (w.r.t. $\mb P_{\theta^\ast}^{(n)}$), 
{\small \begin{align}\label{pf_sketch_1}
\int_{\Theta} \bigg[ \underbrace{\exp \big\{ \ell_n(\theta, \theta^\ast) + n d^{\kappa}(\theta, \theta^\ast) \big\}}_{\xi(\theta, \theta^\ast)} \bigg] \,p_{\theta}(d\theta) \le e^{C n \varepsilon_n^{\kappa}}. 
\end{align}}
To show this, write the integral in \eqref{pf_sketch_1} as $T_1 + T_2$, where $T_1$ and $T_2$ respectively split the integral  over $\{d(\theta, \theta^\ast) \le \varepsilon_n\}$ and $\{d(\theta, \theta^\ast) > \varepsilon_n\}$. Using Markov's inequality along with the fact that $\mb E_{\theta^\ast} [e^{\ell_n(\theta, \theta^\ast)}] = 1$, it can be shown that $T_1 \le e^{C n \varepsilon_n^{\kappa}}$ w.h.p. 
To tackle $T_2$, write $T_2 = T_{21} + T_{22}$ by decomposing $\xi(\theta, \theta^\ast) = \xi(\theta, \theta^\ast) (1 - \phi_n) + \xi(\theta, \theta^\ast) \phi_n$, where $\phi_n$ is the test from Assumption {\bf T}. Using Markov's inequality and the fact that $\phi_n$ has small type-II error, it can be shown that $T_{21} \le e^{-C n \varepsilon_n^{\kappa}}$ w.h.p. The bound on the type-I error of $\phi_n$ along with Markov's inequality yields $\phi_n \le e^{-2 n \varepsilon_n^{\kappa}}$ w.h.p., which yields $T_{22} \le e^{-2 n \varepsilon_n^{\kappa}} T_2$ w.h.p. Combining, one gets $T_2 \le e^{-C n \varepsilon_n^{\kappa}}$ w.h.p. 


Once \eqref{pf_sketch_1} is established, the next step is to link the integrand in \eqref{pf_sketch_1} with the latent variables. To that end, observe that  
$$
\xi(\theta, \theta^\ast) = \sum_{s^n} \exp \{h(\theta, s^n)\} \, \qhat_{S^n}(s^n),
$$
where 
\begin{align*}
h(\theta, s^n) = \log \frac{p(Y^n\,|\,\mu, s^n) \, \pi_{s^n}}{p(Y^n\,|\,\theta^\ast) \, \qhat_{S^n}(s^n)} + n d^{\kappa}(\theta, \theta^\ast). 
\end{align*}
Combining the above with \eqref{pf_sketch_1}, we have, w.h.p.,
\begin{align}\label{eq:pf_sketch_2}
\int_{\Theta} \sum_{s^n} \exp \big\{ h(\theta, s^n) \big\} \, \qhat_{S^n}(s^n) \, p_{\theta}(d \theta) \le e^{C n \varepsilon_n^{\kappa}}. 
\end{align}
Next, use a well-known variational/dual representation of the KL divergence (see, e.g., Corollary 4.15 of \cite{boucheron2013concentration}) which states that for any probability measure $\mu$ and any measurable function $h$ with $e^h \in L_1(\mu)$, 
{\small\begin{align}\label{eq:var_lemma}
\log \int e^{h(\eta)} \, \mu(d \eta) = \sup_{\rho} \bigg[ \int h(\eta) \, \rho(d \eta) -  D(\rho\,\big|\big|\,\mu) \bigg],
\end{align}}
where the supremum is over all probability measures $\rho \ll \mu$. In the present context, setting $\eta = (\theta, s^n)$, $\mu :\,= \qhat_{S^n} \otimes p_{\theta}$, and $\rho = \qhat_{\theta} \otimes \qhat_{S^n}$, it follows from the variational lemma \eqref{eq:var_lemma} and some rearrangement of terms that w.h.p.
{\small \begin{align*}
& n \int_{\Theta} d^{\kappa}(\theta, \theta^\ast) \, \qhat_{\theta}(d \theta) \le n \varepsilon_n^{\kappa} + D(\qhat_{\theta}\,\big|\big|\,p_{\theta})  - \int_{\Theta} \sum_{s^n} h(\theta, s^n) \, \qhat_{\theta}(d \theta).
\end{align*}}
It then follows from \eqref{eq:KL_decomp}--\eqref{eq:V} that the right hand side of the above display equals $n \varepsilon_n^{\kappa} + \Omega(\qhat_{\theta}, \qhat_{S^n})$. The proof of the theorem then follows, since by definition, $\Omega(\qhat_{\theta}, \qhat_{S^n}) \le \Omega(q_{\theta}, q_{S^n})$ for any $(q_{\theta}, q_{S^n})$ in the variational family $\Gamma$.

Now we discuss choices of good variational distributions $q_{\theta}$ and $q_{S^n}$ for minimizing $\Omega(q_{\theta}, q_{S^n})$, the stochastic component of the variational upper bound  in Theorem~\ref{theorem:main}. We make some additional assumptions first on the augmented likelihood and prior in \eqref{eq:aug_post} for the subsequent development. 
First, assume independent priors on $\mu$ and $\pi$ so that $p_{\theta} = p_{\mu} \otimes p_{\pi}$. 
Next, assume $p(Y^n\,|\,\mu, s^n) = \prod_{i=1}^n p(Y_i\,|\,\mu, s_i)$ and $\pi_{s^n} = \prod_{i=1}^n \pi_{s_i}$ splits into independent components. This implies $p(Y^n\,|\,\theta) = \prod_{i=1}^n p(Y_i\,|\,\theta)$, where $p(Y_i\,|\,\theta) = \sum_{s_i} p(Y_i\,|\,\mu, s_i)\,\pi_{s_i}$. For the variational distribution $q_{S^n}$, we additionally assume a {\em mean-field} decomposition $q_{S^n}(s^n) =\prod_{i=1}^n q_{S_i}(S_i)$.

Recall that the objective function decomposes as $\Omega(q_{\theta}, q_{S^n}) = - \int_{\Theta} \ell_n(\theta, \theta^\ast) \, q_{\theta}(d\theta) + \Delta(q_{\theta}, q_{S^n}) + D(q_{\theta}\,\big|\big|\,p_{\theta})$.
The first model-fit term is an averaged (with respect to the variational distribution) log-likelihood ratio which tends to get small as the variational distribution $q_\theta$ places more mass near the true parameter $\theta^\ast$. The second term arising from the approximation of the likelihood function $\ell_n(\theta)$ by $\widetilde{\ell}_n(\theta)$ will become small under some proper choice of $q_{S^n}$, as we will illustrate in the proof of Theorem~\ref{thm:OI_mixture} below. The last regularization or penalty term prevents over-fitting to the data by constricting the KL divergence between the variational solution and the prior. Consequently, a good variational distribution $q_\theta$ should put all its mass into an appropriately small neighborhood around the truth $\theta^\ast$ so that
the first two terms in $V$ become small; on the other hand, the neighborhood has to be large enough so that the last regularization term is not too large.

Motivated by the above discussion, we follow the development of \cite{yang2017variational} by defining two KL neighborhoods around $(\pi^\ast,\,\mu^\ast)$ with radius $(\varepsilon_\pi,\,\varepsilon_\mu)$ as
\begin{align*}
\m B_n(\pi^\ast,\,\varepsilon_\pi) &= \Big\{D(\pi^\ast\, ||\, \pi) \leq \varepsilon_\pi^2,\ \ V(\pi^\ast\, ||\, \pi) \leq \varepsilon_\pi^2\Big\}, \\
\m B_n(\mu^\ast,\,\varepsilon_\mu) &=  \Big\{\sup_{s}D\big[ p(\cdot\,|\,\mu^\ast,s)\, \big|\big| \, p(\cdot\,|\,\mu,s)\big]\leq \varepsilon_\mu^2, \\
&\qquad \sup_s V\big[ p(\cdot\,|\,\mu^\ast,s)\, \big|\big| \, p(\cdot\,|\,\mu,s)\big] \leq \varepsilon_\mu^2\Big\},
\end{align*}
where we used the shorthand $D(\pi^\ast\, ||\, \pi)=\sum_s \pi^\ast_s\,\log(\pi^\ast_s/\pi_s)$ to denote the KL divergence between categorical distributions with parameters $\pi^\ast\in\m S_K$ and $\pi\in\m S_K$ in the $K$-dim simplex $\m S_K$. Consistent with notation introduced at the beginning of Section 2, $V(\pi^\ast\, ||\, \pi)  = \sum_s \pi^\ast_s\,\log^2(\pi^\ast_s/\pi_s)$. 
\begin{theorem}\label{thm:OI_mixture}
For any fixed $(\varepsilon_\pi,\,\varepsilon_\mu) \in (0, 1)^2$, with $\bbP_{\theta^\ast}$ probability at least $1 - 5/\{(D-1)^2 \,n \, (\varepsilon_\pi^2+\varepsilon_\mu^2)\}$, it holds that
\begin{align*}
 & \int d^{\kappa}(\theta, \theta^\ast)\,\qhat_\theta(\theta)\,d\theta \le \varepsilon_n^{\kappa} + (\varepsilon_\pi^2+\varepsilon_\mu^2) + \\
 &\bigg\{ - \frac{\log P_\pi\big[\m B_n(\pi^\ast,\,\varepsilon_\pi)\big]}{n}  \bigg\}  + 
 \bigg\{ - \frac{\log P_\mu\big[B_n(\mu^\ast,\,\varepsilon_\mu) \big] }{n} \bigg\}.
\end{align*}
\end{theorem}
The proof follows Theorem 4.5 of \cite{yang2017variational}, and is omitted; we provide a sketch here. As discussed above, we first make a good choice of $q_{S^n}$ as follows. Let $\qtil_{S^n}$ be a probability distribution over $S^n$ defined as 
\begin{align}\label{Eqn:q_S}
\qtil_{S^n}(s^n) = \prod_{i=1}^n \qtil_{S_i}(s_i)=\prod_{i=1}^n \frac{\pi^\ast_{s_i}\, p(Y_i\,|\,\mu^\ast, \,s_i)}{p(Y_i\,|\,\theta^\ast)}.
\end{align}
Intuitively, $\qtil$ takes the full conditional distribution of $S^n \mid \theta, Y^n$, and replaces $\theta$ by the true parameter $\theta^\ast$. With this choice, $\Delta(q_{\theta}, \qtil_{S^n})$ simplifies to 
\begin{align*}
& - \int_\Theta  q_\theta(\theta)\, \sum_{i=1}^n \sum_{s_i}\qtil_{S_i}(s_i)\,  \log\frac{p(Y_i\,|\,\mu,s_i) \, \pi_{s_i}}{p(Y_i\,|\,\mu^\ast,s_i) \, \pi^\ast_{s_i}}\, d\theta \\
& + \, \int_{\Theta} q_{\theta}(\theta) \,\log\frac{p(Y^n\,|\,\theta)}{p(Y^n\,|\,\theta^\ast)}\,d\theta. 
\end{align*}
It now remains to choose $q_{\theta}$. The first term in the above display naturally suggests choosing 
$q_\theta$ as the restriction of $p_\theta$ into $\m B_n(\pi^\ast,\,\varepsilon_\pi)\times \m B_n(\mu^\ast,\,\varepsilon_\mu)$. For this choice, the second term can also be controlled w.h.p., leading to the conclusion in Theorem \ref{thm:OI_mixture}. 


\section{Applications}
As described in Section 3, variational risk bounds for the parameter of interest depend on the existence of appropriate test functions which characterizes the ability of the likelihood to identify the parameter.  
Developing test functions for studying convergence rates of estimators in classical and Bayesian statistics dates back to \cite{schwartz1965bayes,lecam1973convergence}, with renewed attention in the Bayesian context due to \cite{ghosal2000}. Specific tests have been constructed for nonparametric density estimation \cite{ghosal2000,pati2013posterior}, semi/non-parametric regression \cite{amewou2003posterior,pati2015optimal}, 
convergence of latent mixing measures \cite{nguyen2013convergence,ho2016strong}, high-dimensional problems \cite{pati2014posterior,gao2015rate,chen2016posterior}, and empirical Bayes methods \cite{rousseau2017asymptotic}, among others. As long as the prior distributions are supported on compact subsets of the parameter space, these existing tests can be used to prove minimax optimality of the variational estimate $\hat{\theta}_{VB}$ in each of these problems. We skip the details for space constraints. 

%
In this section, we focus on two novel examples involving latent variables where variational methods are commonly used and no theoretical guarantee is available for the variational solutions. The first one is the Latent Dirichlet Allocation (LDA; \cite{blei2003latent}), a generative probabilistic model for topic modeling. The second example is concerned with estimating the component specific parameters in Gaussian mixture models.


\subsection{Latent Dirichlet allocation}
We first consider LDA \cite{blei2003latent}, a conditionally conjugate probabilistic topic model~\cite{blei2012probabilistic} for learning the latent ``topics'' contained in a collection of documents. Starting from the original paper~\cite{blei2003latent}, the mean-field variational Bayes approximation has become a routine approach for implementing LDA. However, theoretical guarantees for the variational solution is still an open problem despite its empirical success. In this subsection, we show the rate optimality of the estimate from the mean-field approximation to LDA.

In LDA, each document is assumed to contain multiple topics, where a topic is defined as a
distribution over words in a vocabulary.
Our presentation of the model follows the notation of \cite{hoffman2013stochastic}. Let $K$ be the total number of topics, $V$ the vocabulary size, $D$ the total number of documents, and $N$ the number of words in each document (for simplicity, we assume the same number of words across documents). Recall that we use the notation $\mathcal{S}_d$ to denote the $d$-dim simplex. LDA contains two parameters: word distribution matrix $B=[\beta_1,\,\beta_2,\,\ldots,\,\beta_K]\in \mathbb{R}^{V \times K}$ and topic distribution matrix $T_{\gamma} = [\gamma_1,\,\gamma_2,\ldots,\,\gamma_D]\in\mathbb{R}^{K\times D}$, where $\beta_k\in\m S_V$ is the word proportion vector of the $k$th topic, $k=1,2,\dots,K$, and $\gamma_d\in\m S_K$ is the topic proportion vector of the $d$th document, $d=1,2,\ldots,D$.
Given parameters $B$ and $T_{\gamma}$, the data generative model of LDA is:
    \begin{enumerate}
      \item for each document in $d=1,\ldots,D$, draw a topic assignment $z_{dn} \sim  \mbox{Cat}(\gamma_d)$, then
      \item for each word in $n=1,\ldots, N$, draw a word $w_{dn} \sim \mbox{Cat}(\beta_{z_{dn}})$.
  \end{enumerate}
 where $\mbox{Cat}(\pi)$ stands for categorical distribution with probability $\pi$. 
Here, $z_{dn}\in\{1,\ldots,K\}$ is the latent class variable over topics so that
$z_{dn} = k$ indicates the $n$th word in document $d$ is assigned to the
$k$th topic. Similarly, $w_{dn}\in\{1,\ldots,V\}$ is the latent class variable over the words in the vocabulary so that $w_{dn} = v$ indicates
that the $n$th word in document $d$ is the $v$th word in the vocabulary. A common prior distribution over the parameters $B$ and $T_{\gamma}$ is:
  \begin{enumerate}
    \item For each topic in $k=1,\ldots,K$, word proportion vector has prior $\beta_k \sim \mbox{Dir}_V(\eta_\beta)$,
    \item For each document in $d=1,\ldots,D$, topic proportion vector has prior $\gamma_d \sim \mbox{Dir}_K(\eta_\gamma)$.
  \end{enumerate}
Here, $\eta_\beta \in \mb R_{+}$ is a hyper-parameter of the symmetric
Dirichlet prior on the topics $\beta$, and $\eta_\gamma \in \mb R_{+}^K$
are hyper-parameters of the Dirichlet prior on the topic proportions for each document.
To facilitate adaptation to sparsity using Dirichlet distributions when $V,\, K\gg1$, we choose $\eta_\beta=1/V^{c}$ and $\eta_\gamma=1/K^{c}$ for some fixed number $c>1$ \cite{yang2014minimax}.

To apply our theory, we view $N$ as the sample size, and $D$ as the ``dimension" of the parameters in the model.
Under our vanilla notation, we are interested in learning parameters $\theta=(\pi,\,\mu)$, with $\pi=\Gamma$ and $\mu=D$, from the posterior distribution $P(\pi,\, \mu,\,z\,|\, Y^n)$, where $S^{N}=\{S_n:\,n=1,\ldots,N\}$ with $S_n=\{z_{dn}:\, d=1,\ldots, D\}$ are latent variables, $Y^{N}=\{Y_n:\,n=1,\ldots,N\}$ with $Y_n=\{w_{dn}:\, d=1,\ldots, D\}$ are the data, and the priors for $(\pi,\,\mu)$ are independent Dirichlet distributions $\mbox{Dir}_V(\eta_\beta)$ and $\mbox{Dir}_K(\eta_\gamma)$ whose densities are denoted by $p_{\pi}$ and $p_{\mu}$. The conditional distribution $p(Y^N\,|\,\mu,\,S^N)$ of the observation given the latent variable is 
\begin{align*}
\big(w_{dn}\,|\,\mu,\,z_{dn} \big) \sim \mbox{Cat}(\beta_{z_{dn}}),\ d\in[D]\ \mbox{and} \ n\in[N].
\end{align*}
We consider the following mean-field approximation \cite{blei2003latent} by decomposing the variational distribution into
\begin{align*}
  q(\mu, \pi, S^N)
  &=q_\pi(\pi) \, q_\mu(\mu)\, \prod_{n=1}^Nq_{S_n}(S_n)\\
  &=
  \prod_{k=1}^K q_{\beta_k}(\beta_k)
  \prod_{d=1}^D
  \left(
  q_{\theta_d}(\theta_{d})
  \prod_{n=1}^{N}
  q_{z_{dn}}(z_{dn})
  \right)
\end{align*}
for approximating the joint posterior distribution of $(\mu,\,\pi,\, S^N)$. Since for LDA, each observation $Y_n$ is composed of $D$ independent observations, it is natural to present the variational oracle inequality with respect to the average squared Hellinger distance $D^{-1}\sum_{d=1}^D h^2\big[p_d(\cdot\,|\,\theta)\,\big|\big|\,p_d(\cdot\,|\, \theta^\ast)\big]$, where $p_d(\cdot\,|\,\theta)$ denotes the likelihood function of the $d$th observation $w_{d\cdot}$ in $Y_\cdot$. We make the following assumption.

\paragraph{{\bf Assumption S}: (sparsity and regularity condition)} Suppose for each $k$, $\beta^\ast_k$ is $d_k\ll V$ sparse, and for each $d$, $\gamma^\ast_d$ is $e_d\ll K$ sparse. Moreover, there exists some constant $\delta_0>0$, such that each nonzero component of $\beta^\ast_k$ or $\gamma^\ast_d$ is at least $\delta_0$.

\begin{theorem}\label{coro:LDA}
Suppose Assumption {\bf S} holds. If $\big(\sum_{d=1}^De_d+\sum_{k=1}^K d_k\big)/(DN)\to 0$ as $N\to\infty$, then it holds with probability tending to one that as $N\to\infty$
\begin{align*}
& \int D^{-1}\sum_{d=1}^D h^2\big[p_d(\cdot\,|\,\theta)\,\big|\big|\,p_d(\cdot\,|\, \theta^\ast)\big] \ \qhat_\theta(\theta)\,d\theta \\
&\quad \lesssim  \frac{\sum_{d=1}^De_d}{DN}\,\log(DKN) + \frac{\sum_{k=1}^K d_k}{DN}\,\log(KVN).
\end{align*}
\end{theorem}

Theorem~\ref{coro:LDA} implies the estimation consistency as long as the ``effective" dimensionality $\sum_{d=1}^De_d+\sum_{k=1}^K d_k$ of the model is $o(DN)$ as the ``effective sample size" $DN\to\infty$.
In addition, the upper bound depends only logarithmically on the vocabulary size $V$ due to the sparsity assumption.

\subsection{Gaussian mixture models}
Variational inference methods are routinely used in conjugate exponential-family mixtures \cite{ghahramani2001propagation} to speed up computation and perform inference on component-specific parameters, with Gaussian mixtures constituting an important special case. Traditional MCMC methods face difficulties in inferring component-specific parameters due to label-switching \cite{stephens2000dealing}. 
It has been empirically verified \cite{bishop2006pattern} that variational inference for Gaussian mixtures provides accurate estimates of the true density as well as the labels (up to permutation of the indices). However, theoretical guarantees of such a phenomenon is an open problem till date.  In this section, we close this gap and provide an affirmative answer under reasonable assumptions on the true mixture density.  
In particular, we show that variational techniques using mean field approximation provide optimal estimates for the component specific parameters up to some permutation of the labels. 
 
Suppose the true data generating model is the $d$-dimensional Gaussian mixture model with $K$ components,
\begin{eqnarray}
Y_i \mid \mu, \pi \sim \sum_{k=1}^K \pi_k \, \m N(\mu_k,\,  I_d), \ \ i=1, \ldots, n,  \label{eq:gmm}
\end{eqnarray}
where $\mu_k$  is the mean vector associated with the $k$th component and $\pi=(\pi_1,\ldots,\pi_K)\in\m S^{K-1}$ is a vector of probabilities lying in the $K$-dimensional simplex $\m S^{K-1}$.  Assume that $\mu_k\in C_\mu $ where $C_\mu$ is a compact  subset of $\mathbb{R}^d$.  Set $\mu = (\mu_1, \ldots, \mu_K)$ and $\theta = (\mu, \pi)$ as before.  Although we assume the covariance matrix of each Gaussian component to be $I_d$, it is straightforward to extend our results to diagonal covariance matrices. 

Introducing independent latent variables $S_i \in \{1, \ldots, K\}$ for $i=1, \ldots, n$ such that $\pi_k = P(S_i = k)$ for $k=1, \ldots, K$, \eqref{eq:gmm} can be re-written as 
\begin{align}\label{eq:gmml1}
&Y_i \mid \mu, \pi, S^n \sim \m N(\mu_{S_i},\,  I_d),  \\ &S_i \mid \theta \sim \mbox{Multinomial}(1; \pi), \, i=1, \ldots, n, \label{eq:gmml2}
\end{align}

where $S^n = (S_1, \ldots, S_n)$ and $\mbox{Multinomial}(1; \pi)$ denotes a discrete distribution with support $\{1, \ldots, K\}$ and probabilities $\pi_1, \ldots, \pi_K$. For simplicity, we assume independent priors $p_\mu\otimes p_\pi$ for $(\mu,\,\pi)$.

We apply the mean field approximation by using the family of density functions of the form
\begin{align*}
q(\pi, \mu, S^n) = q_\pi(\pi)\, q_\mu(\mu)\, q_{S^n}(s^n) = q_\pi(\pi)\, q_\mu(\mu)\, \prod_{i=1}^n q_{S_i}(s_i)
\end{align*}
to approximate the joint posterior distribution of $(\pi,\,\mu,\, S^n)$.  

Since we are interested in studying accuracy in estimating the component specific parameters $(\pi, \mu)$,  we turn our attention to relating the discrepancy in estimating $f$ with that of $(\pi, \mu)$.  We work with the Wasserstein metric \cite{nguyen2013convergence}  between the mixing measures associated with the density.  Note that \eqref{eq:gmml1}-\eqref{eq:gmml2} can also be written in terms of the mixing measure $P = \sum_{k=1}^K \pi_k \delta_{\mu_k}$ as  
\begin{eqnarray}\label{eq:gmmd}
Y_i \mid P \sim f (\cdot \mid P) = \int \m N( \cdot\,; \mu, \, I_d) \,dP(\mu).  
\end{eqnarray}
Henceforth $\{\mu_k:k = 1, \ldots, K\}$ will be referred to as the atoms of $P$.  Let $\mathcal{P}$ denote the class of all such mixing measures 
\begin{align*}
\m P = \bigg\{P = \sum_{k=1}^K \pi_k \delta_{\mu_k} : \pi \in S^{K-1}, \, \mu_k \in C_{\mu} \, \forall \ k\bigg\},
\end{align*}
such that the atoms lie in a compact subset $C_\mu$ of $\mathbb{R}^d$.  Define the $L_r$-Wasserstein distance, denoted $W_r$, between two mixing measures $P_1 = \sum_{k=1}^K \pi_{1,k} \delta_{\mu_{1,k}}$ and $P_2 = \sum_{k=1}^K \pi_{2,k} \delta_{\mu_{2,k}}$  in $\mathcal{P}$ as 
\begin{eqnarray}\label{eq:W}
W_r^r(P_1, P_2) = \inf_{C \in C_{X_1X_2}} \bbE \|X_1 - X_2\|^r,
\end{eqnarray}
where $X_i \sim P_i$ for $i = 1, 2$, and $C_{X_1 X_2}$ is the set of all possible {\em couplings}, i.e. joint distributions of $X_1$ and $X_2$ with marginals $P_1$ and $P_2$ respectively.  One can write \eqref{eq:W} in terms of $(\pi, \mu)$ as 
\begin{eqnarray}\label{eq:W2}
W_r^r(P_1, P_2) = \inf_{q \in C_{P_1P_2}} \sum_{k, k'} q_{kk'} \|\mu_{1,k} -  \mu_{2,k'}\|^r,
\end{eqnarray}
where $q$ varies over $C_{P_1P_2}$, the set of joint probability mass functions over $\{1, \ldots, K\}^2$ satisfying $\sum_{k} q_{kk'} = \pi_{2,k'}$ and $\sum_{k'} q_{kk'} = \pi_{1,k}$.   

In the following, we will consider $r=1$ to work with the $W_1$ metric.  It is known \cite{nguyen2013convergence} that $\mathcal{P}$ is compact with respect to $W_1$.  Let $h\{ f(\cdot \mid P_1) \,\big|\big|\, f(\cdot \mid P_2)\}$ denote the Hellinger distance between the densities with corresponding mixing measures $P_1$ and $P_2$.  Denote by $N(\varepsilon, \mathcal{P}, W_1)$ and $N(\varepsilon, \mathcal{F},h)$ the $\varepsilon$-covering numbers of $\mathcal{P}$ and the corresponding space of densities $\mathcal{F}$ with respect to $W_1$ and $h$ respectively.

Following \cite{nguyen2013convergence}, we investigate the minimum separation between the densities in \eqref{eq:gmmd} in Hellinger distance when the corresponding mixing measures are separated by at least $\varepsilon$ in $W_1$.   Fix   $P^\ast = \sum_{k=1}^K \pi_k^\ast \delta_{\mu_k^\ast} \in \mathcal{P}$ and $P \in \mathcal{P}$. 
From Birg\'{e}-Lecam theory \cite{birge1983approximation,lecam1973convergence}, there exists a sequence test functions $\Phi_n$ based on observations $Y_1, \ldots, Y_n$ such that 
\begin{align}\label{eq:BLtest1}
& \bbE_{P^\ast} \Phi_n \leq N(\varepsilon, \mathcal{F}, h) \,e^{-C_1n \varepsilon^2}\\
& \bbE_{P} [ 1- \Phi_n]  \leq e^{-C_2n h^2[f(\cdot\,|\,P)\,||\,f(\cdot\,|\, P^\ast)]},\label{eq:BLtest2}
\end{align} 
 for any $P$ with $h[f(\cdot\,|\,P)\,||\,f(\cdot\,|\, P^\ast)] > \varepsilon$. \cite{ho2016strong} discusses construction of test functions in the $W_1$ metric. Using their Theorem 3.1 along with \eqref{eq:BLtest1}-\eqref{eq:BLtest2}, one obtains a test function $\Psi_n$ such that 
\begin{align}
&\bbE_{P^\ast} [\Psi_{n}] \leq  \, N(\varepsilon,\m P, W_1)\, e^{-C_1 n \varepsilon^2} \label{eq:gmtest1}\\
&\bbE_{P} [ 1- \Psi_{n}]  \leq e^{-C_2 n W_1^2(P, P^\ast)},  \label{eq:gmtest2}
\end{align}
for any $P \in \m P$ with $W_1(P, P^\ast) > \varepsilon$. We show in the supplement that $N(\varepsilon,\m P, W_1) \lesssim e^{C n \varepsilon^2}$. Hence Assumption {\bf T} is satisfied with $\kappa = 2$. We remark that for the $W_2$ metric, it is possible to construct a test with $\kappa = 4$. 

To state the risk bounds in the $W_1$ metric, we require the following assumptions on the component specific weights  associated with the true mixing measure $P^\ast$ and the prior densities $p_\mu$ and $p_\pi$ respectively. 
\paragraph{{\bf Assumption R}: (lower bound on component weights)}  There exists some constant $\delta \in (0, 1)$, such that the true $\pi_k^\ast > \delta$ for all $k = 1, \ldots, K$.

\paragraph{{\bf Assumption P}: (Prior thickness)}  Assume that the prior densities for $\mu$ and $\pi$ satisfy 
$p_\pi(\pi^\ast)>0$ and $p_\mu(\mu^\ast)>0$. Additionally, assume that $p_{\mu}$ is supported on the compact set $C_{\mu}$.

\begin{theorem}\label{thm:GM}
Suppose Assumption {\bf R} holds, and the prior densities $p_\mu$ and $p_\pi$ satisfy Assumption {\bf P}.  Then if  $d\,K/n\to 0$ as $n\to \infty$,  it holds with probability tending to one as $n\to\infty$ that
\begin{eqnarray} \label{eq:1}
 \int h^2\big[p(\cdot\,|\,\theta)\,\big|\big|\, p(\cdot\,|\, \theta^\ast)\big] \  \qhat_\theta(\theta)\,d\theta  \lesssim \frac{d\,K}{n}\,\log (dn), \\
 \int W_1^2(P, P^\ast) \ \qhat_\theta(\theta)\,d\theta \lesssim \frac{d\,K}{n}\,\log (K\,n). \label{eq:3}
\end{eqnarray}
\end{theorem} 
The convergence rates obtained in \eqref{eq:1}-\eqref{eq:3} are minimax upto logarithmic terms  \cite{ho2016strong}.
An important consequence of Theorem \ref{thm:GM} is that it allows us to study the accuracy of the estimates of the component specific parameters $(\pi_k, \mu_k)$ for $k=1, \ldots, K$. The following lemma relates the accuracy in estimating the mixing measures with respect to the $W_1$ distance and that for the component specific weights and atoms.  The convergence of the weights requires an additional assumption in terms of separability of atoms in $P^\ast$ which prohibits two distinct atoms of the true mixing distribution $\m P^\ast$ to be vanishingly close.   

\paragraph{{\bf Assumption S}: (separability of atoms)}  There exists some $\zeta>0$, such that  $\inf_{k\neq k'} \| \mu_k^\ast -   \mu_{k'}^\ast \| \geq \zeta$.

\begin{lemma}\label{lem:cs}
 Fix $P^\ast \in \mathcal{P}$ such that $P^\ast$ satisfies Assumption {\bf R}.  Given $\varepsilon > 0$, if $P = \sum_{k=1}^K \pi_k \delta_{\mu_k}$ satisfies $W_1(P^\ast, P) < \varepsilon$, then  for $\varepsilon < \zeta \delta$, 
\begin{align*}
& \|\mu_k - \mu_k^\ast \| < \varepsilon/\delta, \ \forall \ k=1, \ldots, K,
\end{align*}
upto a permutation of the indices.  
If in addition, $P^\ast$ satisfies  both Assumption {\bf R} and {\bf S}, then,  
\begin{eqnarray*}
\sum_{k=1}^K |\pi_k - \pi^\ast_k| <  \frac{\varepsilon}{ (\zeta - \varepsilon/\delta)}. 
\end{eqnarray*}

\end{lemma}

Theorem \ref{thm:GM} together with Lemma \ref{lem:cs} implies the following corollary about parameter convergence rates.
 \begin{corollary}
 Let  $M_n \uparrow \infty$ be any sequence of positive numbers and  $\varepsilon_n =  \sqrt{(dK /n)\log (K\,n)}$. If 
 $P^\ast$ satisfies Assumption {\bf R},   
 then $ \qhat_\theta \big( \|\mu_{k}^\ast -  \mu_{k}\| < M_n \varepsilon_n, \ \forall \, k\big) \to 1\, \mbox{a.s.}$ 
 If $P^\ast$ satisfies both assumptions {\bf R} and {\bf S}, then  
$
 \qhat_\theta \big( \|\mu_{k}^\ast -  \mu_{k}\| < M_n \varepsilon_n, \ \forall \, k, \,  \sum_{k=1}^K |\pi_k^\ast - \pi_k| < M_n\varepsilon_n \big) \to 1 \, \mbox{a.s.} 
$
 \end{corollary}
 
According to \cite{chen1995}, $n^{-1/2}$ is the minimax rate of estimating $\mu_k$'s when the number of components are correctly specified to be $K$ ($\pi_k$'s are not estimable if some $\mu_k^\ast$'s are the same). The second part of this corollary implies $n^{-1/2}$ convergence rates (up to $\log n$ factors) for both $\mu$ and $\pi$ under the strong identifiability assumption {\bf S}.

\section{Conclusion}

We have provided general purpose tools to verify the statistical accuracy of point estimates obtained from mean-field variational inference. Our analysis incorporates latent variables commonly augmented to simplify variational inference, and hence our theory is applicable to a broad range of existing algorithms in exactly the way they are practically implemented. The theory does not require prior conjugacy or analytical tractability of the iterative algorithms. The only two ingredients are the existence of appropriate tests which quantify the ability of the likelihood to separate points in the parameter space relative to the distance metric $d$, and the prior mass assigned to appropriate KL neighborhoods of the true parameter. Since similar quantities also appear in quantifying the convergence rate of the actual posterior \cite{rousseau2016frequentist}, one can leverage on this large body of literature along with our theory to offer insights into variational point estimates. Future work will involve relaxing the prior compactness assumption and considering more general variational families.  

\section{Proofs} \label{appendix}
\subsection{Convention}
Equations in the main document are cited as (1), (20 etc., retaining their numbers, while new equations defined in this document are numbered (S1), (S2) etc. 

\subsection{Proof of Theorem \ref{theorem:main}}
As in the proof sketch in the main document, our first step is to show that under the testing assumption {\bf T}, 
\begin{align}\label{eq:cond1}
\int_{\Theta} \xi(\theta, \theta^\ast)\,p_{\theta}(d\theta) \le e^{C n \varepsilon_n^\kappa},
\end{align}
w.h.p. (w.r.t. $\mb P_{\theta^\ast}^{(n)}$), where recall $\log \xi(\theta, \theta^\ast) = \ell_n(\theta, \theta^\ast) + n d^{\kappa}(\theta, \theta^\ast)$. We first establish \eqref{eq:cond1}. Define 
\begin{align*}
& T_1 = \int_{d(\theta, \theta^\ast) \le \varepsilon_n} \xi(\theta, \theta^\ast)\,p_{\theta}(d\theta), \\
& T_2 = \int_{d(\theta, \theta^\ast) > \varepsilon_n} \xi(\theta, \theta^\ast)\,p_{\theta}(d\theta). 
\end{align*}
Let us first tackle $T_1$. Since $\mb E_{\theta^\ast} [e^{ \ell_n(\theta, \theta^\ast}] = 1$, we have, 
$$
\mb E_{\theta^\ast} T_1 = \int_{d(\theta, \theta^\ast) \le \varepsilon_n} e^{n d^{\kappa}(\theta, \theta^\ast)} \,p_{\theta}(d\theta) \le e^{n \varepsilon_n^{\kappa}}. 
$$
Hence, by Markov's inequality, $T_1 \le e^{C n \varepsilon_n^{\kappa}}$ with probability at least $1 - e^{- C n \varepsilon_n^{\kappa}}$. 

Let us now focus on $T_2$. Write $T_2 = T_{21} + T_{22}$, where 
\begin{align*}
& T_{21} = \int_{d(\theta, \theta^\ast) > \varepsilon_n} (1 - \phi_n) \,\xi(\theta, \theta^\ast)\,p_{\theta}(d\theta), \\
& T_{22} = \int_{d(\theta, \theta^\ast) > \varepsilon_n} \phi_n \, \xi(\theta, \theta^\ast)\,p_{\theta}(d\theta),
\end{align*}
where $\phi_n$ is the test function from Assumption {\bf T}. Focus on $T_{21}$ first. Observe 
\begin{align*}
\mb E_{\theta^\ast} T_{21} & = \int_{d(\theta, \theta^\ast) > \varepsilon_n} \mb E_{\theta}[1 - \phi_n] \, e^{n d^\kappa(\theta, \theta^\ast)}\,p_{\theta}(d \theta) \\& \le e^{- C n \varepsilon_n^\kappa}.
\end{align*}
This implies, by Markov's inequality, than $T_{21} \le e^{- C n \varepsilon_n^\kappa}$ with probability at least $1 - e^{- C n \varepsilon_n^{\kappa}}$. 

Finally, focus on $T_{22}$. Since $\mb E_{\theta^\ast} [\phi_n] \le e^{-  n \varepsilon_n^{\kappa}}$, it follows from Markov's inequality that $\phi_n \le e^{- C n \varepsilon_n^\kappa}$ with probability at least $1 - e^{- C n \varepsilon_n^{\kappa}}$. Hence, $T_{22} \le e^{-C n \varepsilon_n^\kappa} T_2$ w.h.p. Adding the w.h.p. bound for $T_{21}$, we obtain, w.h.p.,
$$
T_{2} \le e^{-C n \varepsilon_n^\kappa} T_2 + e^{-C n \varepsilon_n^\kappa}. 
$$
Rearranging, $T_2 \le e^{- C n \varepsilon^\kappa}$ with probability at least $1 - e^{- C n \varepsilon_n^{\kappa}}$. Combining with the bound for $T_1$, \eqref{eq:cond1} is established. 

Once \eqref{eq:cond1} is established, the next step is to link the integrand in \eqref{eq:cond1} with the latent variables. To that end, observe that  
$$
\xi(\theta, \theta^\ast) = \sum_{s^n} \exp \{h(\theta, s^n)\} \, \qhat_{S^n}(s^n),
$$
where 
\begin{align*}
h(\theta, s^n) = \log \frac{p(Y^n\,|\,\mu, s^n) \, \pi_{s^n}}{p(Y^n\,|\,\theta^\ast) \, \qhat_{S^n}(s^n)} + n d^{\kappa}(\theta, \theta^\ast). 
\end{align*}
Combining the above with \eqref{pf_sketch_1}, we have, w.h.p.,
\begin{align}\label{eq:pf_sketch_2}
\int_{\Theta} \sum_{s^n} \exp \big\{ h(\theta, s^n) \big\} \, \qhat_{S^n}(s^n) \, p_{\theta}(d \theta) \le e^{C n \varepsilon_n^{\kappa}}. 
\end{align}
Next, use a well-known variational/dual representation of the KL divergence (see, e.g., Corollary 4.15 of \cite{boucheron2013concentration}) which states that for any probability measure $\mu$ and any measurable function $h$ with $e^h \in L_1(\mu)$, 
{\small\begin{align}\label{eq:var_lemma}
\log \int e^{h(\eta)} \, \mu(d \eta) = \sup_{\rho} \bigg[ \int h(\eta) \, \rho(d \eta) -  D(\rho\,\big|\big|\,\mu) \bigg],
\end{align}}
where the supremum is over all probability measures $\rho \ll \mu$. In the present context, setting $\eta = (\theta, s^n)$, $\mu :\,= \qhat_{S^n} \otimes p_{\theta}$, and $\rho = \qhat_{\theta} \otimes \qhat_{S^n}$, it follows from the variational lemma \eqref{eq:var_lemma} and some rearrangement of terms that w.h.p.
{\small \begin{align*}
& n \int_{\Theta} d^{\kappa}(\theta, \theta^\ast) \, \qhat_{\theta}(d \theta) \le n \varepsilon_n^{\kappa} + D(\qhat_{\theta}\,\big|\big|\,p_{\theta})  - \int_{\Theta} \sum_{s^n} h(\theta, s^n) \, \qhat_{\theta}(d \theta).
\end{align*}}
From \eqref{eq:KL_decomp}--\eqref{eq:V} (in the main document), it follows that the right hand side of the above display equals $n \varepsilon_n^{\kappa} + \Omega(\qhat_{\theta}, \qhat_{S^n})$. The proof of the theorem then follows, since by definition, $\Omega(\qhat_{\theta}, \qhat_{S^n}) \le \Omega(q_{\theta}, q_{S^n})$ for any $(q_{\theta}, q_{S^n})$ in the variational family $\Gamma$.

\subsection{Proof of Lemma \ref{lem:cs}}
Since $W_1(P^\ast, P) < \varepsilon$,  there exists a coupling $q$ such that  $ \sum_{k, k'} q_{kk'} \|\mu_{k}^\ast -  \mu_{k'}\| < \varepsilon$.   Then $ \sum_{k} \pi_k^\ast \inf_{k'}\|\mu_{k}^\ast -  \mu_{k'}\| < \varepsilon$. Since $\pi_k^\ast \geq \delta$, we have  $\inf_{k'}\|\mu_{k}^\ast -  \mu_{k'}\| \leq \varepsilon/\delta$ for all $k= 1, \ldots, K$.  
 This means for any $k$, there exists a $k'$ such that $\|\mu_{k}^\ast -  \mu_{k'}\| < \varepsilon/\delta$.  Without loss of generality, let $k' = k$.   This proves the first part of the assertion.   To prove the second part, observe that for $k \neq k'$, $\|\mu_{k}^\ast -  \mu_{k'}\| \geq \zeta -  \|\mu_{k'}^\ast -  \mu_{k'}\| \geq \kappa - \varepsilon/\delta$.  Then 
 \begin{align*}
 &&\varepsilon > W_1(P^\ast, P) \geq  \inf_{q} \sum_{k\neq k'} q_{kk'} \|\mu_{k}^\ast -  \mu_{k'}\| \\ && \geq (\zeta - \varepsilon/\delta)
 \inf_{C \in C_{XY}} \bbP(X \neq Y) \\
 && = (\zeta - \varepsilon/\delta) \sum_{k=1}^K |\pi_k^\ast - \pi_k|,  
 \end{align*}
 implying $\sum_{k=1}^K |\pi_k^\ast - \pi_k| \leq \varepsilon / (\zeta - \varepsilon/\delta)$.

\subsection{Proof of Theorem \ref{thm:GM}}
We first ensure the existence of the test functions $\Phi_n, $ and $\Psi_{n}$ as described in 
\eqref{eq:BLtest1}-\eqref{eq:gmtest2}. 
First, we find find the covering numbers $N(\varepsilon, \mathcal{P}, W_1)$ and $N(\varepsilon, \mathcal{F}, h)$ to upper bound the Type I and II errors of the test functions $\Phi_n$ and $\Psi_n$.  Note that 
\begin{eqnarray*}
 h^2[f(\cdot\,|\,P_1)\,||\,f(\cdot\,|\, P_2)] \leq \sum_{k=1}^K |\pi_{1,k} - \pi_{2,k}| + \\ \sum_{k=1}^k \pi_{1,k} \norm{\mu_{1,k} - \mu_{2,k}}. 
 \end{eqnarray*}
Hencd $N(\varepsilon, \mathcal{F}, h) \leq N(\varepsilon^2/2, \mathcal{S}^{K-1}, || \cdot ||_1) \times \{N(\varepsilon^2/2, C_{\mu}, || \cdot ||)\}^K$ where $|| \cdot ||_1$ denotes the $L_1$ norm between
 two probability vectors and  $|| \cdot ||$ denotes the Euclidean norm.  From Lemma A.4 of \cite{ghosal2001entropies}, we obtain $N(\varepsilon^2/2, \mathcal{S}^{K-1}, || \cdot ||_1) \leq (10/\varepsilon^2)^{K-1}$.  Also,  $\{N(\varepsilon^2/2, C_{\mu}, || \cdot ||)\} \leq  (2C_U/ \varepsilon^2)^d$ for a global constant $C_U$ is the diameter of the set $C_\mu$. Then 
$ N(\varepsilon, \mathcal{F}, h) \leq (C/\varepsilon^2)^{dK}$ for some constant $C > 0$.  
To obtain an upper bound for $N(\varepsilon, \mathcal{P}, W_1)$, we note that 
\begin{eqnarray*}
W_1(P_1, P_2) \leq \sum_{k=1}^K \max\{\pi_{1,k}, \pi_{2,k}\}\norm{\mu_{1,k} - \mu_{2,k}} \\ + C_U\sum_{k=1}^K |\pi_{1,k} - \pi_{2,k}|.
\end{eqnarray*}
Hence $N(\varepsilon, \mathcal{P}, W_1) \leq N(\varepsilon/(2 C_U), \mathcal{S}^{K-1}, || \cdot ||_1) \times \{N(\varepsilon/(2K), C_{\mu}, || \cdot ||)\}^K \leq( CK/\varepsilon)^{dK} (10/\varepsilon)^{K-1}$. 
Hence   $\log N(\varepsilon, \mathcal{F}, h) \lesssim dK \log (1/\varepsilon)$ and $\log N(\varepsilon, \mathcal{P}, W_1) \lesssim dK \log (K/\varepsilon)$. 
Then, we have from \eqref{eq:BLtest1}-\eqref{eq:BLtest2}
\begin{eqnarray}
\bbE_{P^\ast} \Phi_n \leq e^{-C_1n \varepsilon^2 + dK \log (1/\varepsilon)}\\
 \bbE_{P} [ 1- \Phi_n]  \leq e^{-C_2n h^2[f(\cdot\,|\,P)\,||\,f(\cdot\,|\, P^\ast)]},
\end{eqnarray} 
for any $P$ with $h[f(\cdot\,|\,P)\,||\,f(\cdot\,|\, P^\ast)] > \varepsilon$. 
In this case, we choose $\varepsilon \equiv \varepsilon_n$ to be as constant multiple of $\{(dK/n) \log n\}^{1/2}$.  
Also, we  have from  \eqref{eq:gmtest1}--\eqref{eq:gmtest2}
\begin{eqnarray}
\bbE_{P^\ast} \Psi_n \leq e^{-C_1 n \varepsilon^2 +  dK \log (K/\varepsilon)}\\
 \bbE_{P} [ 1- \Psi_n]  \leq e^{-C_2 n W_1^2(P, P^\ast)},
\end{eqnarray}
for any $P$ with $W_1(P, P^\ast) > \varepsilon$. 
In this case, we choose $\varepsilon \equiv \varepsilon_n$ to be as constant multiple of $\{(dK/n) \log (Kn)\}^{1/2}$.

Recall the two KL neighborhoods around $(\pi^\ast,\,\mu^\ast)$ with radius $(\varepsilon_\pi,\,\varepsilon_\mu)$ as
\begin{align*}
\m B_n(\pi^\ast,\,\varepsilon_\pi) &= \Big\{D(\pi^\ast\, ||\, \pi) \leq \varepsilon_\pi^2,\ \ V(\pi^\ast\, ||\, \pi) \leq \varepsilon_\pi^2\Big\}, \\
\m B_n(\mu^\ast,\,\varepsilon_\mu) &=  \Big\{\sup_{s}D\big[ p(\cdot\,|\,\mu^\ast,s)\, \big|\big| \, p(\cdot\,|\,\mu,s)\big]\leq \varepsilon_\mu^2, \\
&\qquad \sup_s V\big[ p(\cdot\,|\,\mu^\ast,s)\, \big|\big| \, p(\cdot\,|\,\mu,s)\big] \leq \varepsilon_\mu^2\Big\},
\end{align*}
where we used the shorthand $D(\pi^\ast\, ||\, \pi)=\sum_s \pi^\ast_s\,\log(\pi^\ast_s/\pi_s)$ to denote the KL divergence between multinomial distributions with parameters $\pi^\ast, \pi \in\m S_K$. We choose $q_\theta$ as the restriction of $p_\theta$ into $\m B_n(\pi^\ast,\,\varepsilon_\pi)\times \m B_n(\mu^\ast,\,\varepsilon_\mu)$. 


It is easy to verify that under Assumption {\bf R}, there exists some constant $C_1$ depending only on $\delta_0$ such that 
$\m B_n(\pi^\ast,\,\sqrt{K}\,\varepsilon)\supset \{\pi:\,\max_{k}|\pi_k-\pi^\ast_k|\leq C_1\,\varepsilon\}$ (by using the inequality $D(p\,||\,q) \geq 2\,h^2(p\,||\,q)$). In addition, for Gaussian mixture model, it is easy to verify that the KL neighborhood $\m B_n(\mu^\ast,\,\varepsilon)$  contains the set $\{\mu:\,\max_{k}\|\mu_k-\mu_k^\ast\|\leq 2\,\varepsilon\}$. As a consequence,  with $\varepsilon_\pi=\sqrt{K}\,\varepsilon$ and $\varepsilon_\mu=\varepsilon$ yields (using the prior thickness assumption and the fact that the volumes of $\{\pi:\,\max_{k}|\pi_k-\pi^\ast_k|\leq C_1\,\varepsilon\}$ and $\{\mu:\,\max\|\mu_k-\mu^\ast_k\|\leq C_2\,\varepsilon\}$ are at least $\m O(\varepsilon^{-K})$ and $\m O\big((\sqrt{d}/\varepsilon)^{dK}\big)$ respectively). Then we have from   Theorem \ref{thm:OI_mixture},with probability tending to one as $n\to\infty$,
\begin{align*}
 \int \Big\{h^2\big[f(\cdot\,|\,\theta)\,\big|\big|\,f(\cdot\,|\, \theta^\ast)\big] \Big\}\,\qhat_\theta(\theta)\,d\theta &\lesssim \frac{d\,K}{n}\log n+K\,\varepsilon^2  \\ & + \frac{d\,K}{n} \,\log\frac{d}{\varepsilon}.
\end{align*}
Choosing $\varepsilon = \sqrt{d/n}$ in the above display yields the claimed bound.

Also, we have with high probability 
\begin{align*}
 \int \Big\{W_1^2\big[f(\cdot\,|\,\theta)\,\big|\big|\,f(\cdot\,|\, \theta^\ast)\big] \Big\}\,\qhat_\theta(\theta)\,d\theta &\lesssim \frac{d\,K}{n}\log (Kn) \\ & + K\,\varepsilon^2+ \frac{d\,K}{n} \,\log\frac{d}{\varepsilon}.
\end{align*}
Choosing $\varepsilon = \sqrt{d/n}$ in the above display yields the claimed bound noting that the first term in the right hand side of the preceding display is dominant.

\subsection{Proof of Theorem~\ref{coro:LDA}}
Under the notation in the paper, for each $n=1,\ldots,N$, the latent variable $S_n=\{z_{dn}:\, d=1,\ldots,D\}$. We use Theorem~\ref{thm:OI_mixture} with $d=h$ (Hellinger metric) and view each latent variable $S_n$ per observation in the theorem as a block of $D$ independent latent variable per observation. The existence of the test is automatic \cite{lecam1973convergence} with the Hellinger metric (parameter space is compact). This leads to that with probability tending to one as $N\to \infty$,
\begin{align*}
 &\int \sum_{d=1}^D h^2\big[p_d(\cdot\,|\,\theta)\,\big|\big|\,p_d(\cdot\,|\, \theta^\ast)]\,d\theta
  \le \bigg(\sum_{d=1}^D\varepsilon_{\gamma_d}^2+\sum_{k=1}^K\varepsilon_{\beta_k}^2\bigg)\\
 &\qquad + \Big\{ - \frac{1}{N} \sum_{d=1}^D\log P_{\gamma_d}\big[\m B_N(\gamma^\ast_d,\,\varepsilon_{\gamma_d})\big] \Big\} \\
 &\qquad+ \Big\{ - \frac{1}{N} \sum_{k=1}^K\log P_{\beta_k}\big[B_N(\beta^\ast_k,\,\varepsilon_{\beta_k}) \big] \Big\},
\end{align*}
where KL neighborhoods $\m B_N(\gamma_d^\ast;\, \varepsilon_{\gamma_d}):\,=\big\{D(\gamma_d^\ast\,||\,\gamma_d)$ $\leq \varepsilon_{\gamma_d}^2,$ $V(\gamma_d^\ast\,||\,\gamma_d)\leq \varepsilon_{\gamma_d}^2\big\}$, for $d=1,\ldots,D$, and  $B_N(\beta^\ast_k,\,\varepsilon_{\beta_k}) = \big\{\max_k D\big[p(\cdot\,|\, \beta_k,\,k)\,||\, p(\cdot\,|\, \beta_k,\,k)\big]\leq \varepsilon_{\beta_k}^2,\, \max_{S_n} V\big[p(\cdot\,|\, \beta_k,\,k)\,||\, p(\cdot\,|\, \beta_k,\,k)\big]\leq \varepsilon_{\beta_k}^2\big\}$.

Let $S^\beta_k$ denote the index set corresponding to the non-zero components of $\beta_k$ for $k=1,\ldots,K$, and $S^\gamma_d$ the index set corresponding to the non-zero components of $\gamma_d$ for $d=1,\ldots,D$.
Under Assumption {\bf S}, it is easy to verify that for some sufficiently small constants $c_1,c_2>0$, it holds for all $d=1,\ldots,D$ that $\m B_N(\gamma^\ast_d,\,\varepsilon_{\gamma_d})\supset \big\{ \|(\gamma_d)_{(S^\gamma_d)^c}\|_1 \leq c_1\, \varepsilon_{\gamma_d}, \, \|(\gamma_d)_{S^\gamma_d}- (\gamma^\ast_d)_{S^\gamma_d}\|_{\infty} \leq c_1\,\varepsilon_{\gamma_d}\big\}$, and for all $k=1,\ldots,K$ that $\m B_N(\beta^\ast_k,\,\varepsilon_{\beta_k})\supset \big\{ \|(\beta_k)_{(S^\beta_k)^c}\|_1 \leq c_2\, \varepsilon_{\beta_k}, \, \|(\beta_k)_{S^\beta_k}- (\beta^\ast_k)_{S^\beta_d}\|_{\infty} \leq c_2\,\varepsilon_{\beta_d}\big\}$.
Applying Theorem 2.1 in \cite{yang2014minimax}, we obtain the following prior concentration bounds for high-dimensional Dirichlet priors 
\begin{align*}
&P_{\gamma_d}\big\{ \|(\gamma_d)_{(S^\gamma_d)^c}\|_1 \leq c_1\, \varepsilon_{\gamma_d}, \, \|(\gamma_d)_{S^\gamma_d}- (\gamma^\ast_d)_{S^\gamma_d}\|_{\infty} \leq c_1\,\varepsilon_{\gamma_d}\big\} \\
&\qquad\qquad \qquad \gtrsim \exp\Big\{-C\,e_d\,\log\frac{K}{\varepsilon_{\gamma_d}} \Big\},\ d=1,\ldots,D;\\
&P_{\beta_k}\big\{ \|(\beta_k)_{(S^\beta_k)^c}\|_1 \leq c_2\, \varepsilon_{\beta_k}, \, \|(\beta_k)_{S^\beta_k}- (\beta^\ast_k)_{S^\beta_k}\|_{\infty} \leq c_2\,\varepsilon_{\beta_d}\big\}\\ &\qquad\qquad \qquad\gtrsim \exp\Big\{-C\,d_k\,\log\frac{V}{\varepsilon_{\beta_k}} \Big\}, \ k=1,\ldots,K,
\end{align*}
for some constant $C>0$.

Putting pieces together, we obtain
\begin{align*}
 &\int \sum_{d=1}^D h^2\big[p_d(\cdot\,|\,\theta)\,\big|\big|\,p_d(\cdot\,|\, \theta^\ast)]\,d\theta
 \lesssim  \bigg(\sum_{d=1}^D\varepsilon_{\gamma_d}^2+\sum_{k=1}^K\varepsilon_{\beta_k}^2\bigg) \\
 &+ \frac{1}{N} \sum_{d=1}^D e_d\,\log\frac{K}{\varepsilon_{\gamma_d}} +
  \frac{1}{N} \sum_{k=1}^K d_k\,\log\frac{V}{\varepsilon_{\beta_k}},
\end{align*}
which leads to the desired bound by optimally choosing $\varepsilon_{\gamma_d}$'s and $\varepsilon_{\beta_k}$'s.

\bibliographystyle{plain}

\bibliography{VB}

\begin{thebibliography}{10}

\bibitem{ahmed2012scalable}
Amr Ahmed, Mohamed Aly, Joseph Gonzalez, Shravan Narayanamurthy, and Alexander
  Smola.
\newblock Scalable inference in latent variable models.
\newblock In {\em International conference on Web search and data mining
  (WSDM)}, volume~51, pages 1257--1264, 2012.

\bibitem{alquier2017concentration}
Pierre Alquier and James Ridgway.
\newblock Concentration of tempered posteriors and of their variational
  approximations.
\newblock {\em arXiv preprint arXiv:1706.09293}, 2017.

\bibitem{alquier2016properties}
Pierre Alquier, James Ridgway, and Nicolas Chopin.
\newblock On the properties of variational approximations of gibbs posteriors.
\newblock {\em JMLR}, 17(239):1--41, 2016.

\bibitem{amewou2003posterior}
Messan Amewou-Atisso, Subhashis Ghosal, Jayanta~K Ghosh, and RV~Ramamoorthi.
\newblock Posterior consistency for semi-parametric regression problems.
\newblock {\em Bernoulli}, 9(2):291--312, 2003.

\bibitem{attias2000variational}
Hagai Attias.
\newblock {A variational Bayesian framework for graphical models}.
\newblock In {\em Advances in neural information processing systems}, pages
  209--215, 2000.

\bibitem{birge1983approximation}
Lucien Birg{\'e}.
\newblock Approximation dans les espaces m{\'e}triques et th{\'e}orie de
  l'estimation.
\newblock {\em Probability Theory and Related Fields}, 65(2):181--237, 1983.

\bibitem{bishop2006pattern}
Christopher~M Bishop.
\newblock {\em Pattern recognition and machine learning}.
\newblock springer, 2006.

\bibitem{blei2012probabilistic}
David~M Blei.
\newblock Probabilistic topic models.
\newblock {\em Communications of the ACM}, 55(4):77--84, 2012.

\bibitem{blei2017variational}
David~M Blei, Alp Kucukelbir, and Jon~D McAuliffe.
\newblock Variational inference: A review for statisticians.
\newblock {\em Journal of the American Statistical Association},
  (just-accepted), 2017.

\bibitem{blei2003latent}
David~M Blei, Andrew~Y Ng, and Michael~I Jordan.
\newblock Latent dirichlet allocation.
\newblock {\em Journal of machine Learning research}, 3(Jan):993--1022, 2003.

\bibitem{boucheron2013concentration}
St{\'e}phane Boucheron, G{\'a}bor Lugosi, and Pascal Massart.
\newblock {\em Concentration inequalities: A nonasymptotic theory of
  independence}.
\newblock Oxford university press, 2013.

\bibitem{casella2002statistical}
George Casella and Roger~L Berger.
\newblock {\em Statistical inference}, volume~2.
\newblock Duxbury Pacific Grove, CA.

\bibitem{chen1995}
Jiahua Chen.
\newblock Optimal rate of convergence for finite mixture models.
\newblock {\em Ann. Statist.}, 23(1):221--233, 1995.

\bibitem{chen2016posterior}
Mengjie Chen, Chao Gao, Hongyu Zhao, et~al.
\newblock {Posterior Contraction Rates of the Phylogenetic Indian Buffet
  Processes}.
\newblock {\em Bayesian analysis}, 11(2):477--497, 2016.

\bibitem{corduneanu2001variational}
Adrian Corduneanu and Christopher~M Bishop.
\newblock {Variational Bayesian model selection for mixture distributions}.
\newblock In {\em Artificial intelligence and Statistics}, volume 2001, pages
  27--34. Waltham, MA: Morgan Kaufmann, 2001.

\bibitem{gao2015rate}
Chao Gao, Harrison~H Zhou, et~al.
\newblock Rate-optimal posterior contraction for sparse pca.
\newblock {\em The Annals of Statistics}, 43(2):785--818, 2015.

\bibitem{gelfand1990sampling}
Alan~E Gelfand and Adrian~FM Smith.
\newblock Sampling-based approaches to calculating marginal densities.
\newblock {\em Journal of the American statistical association},
  85(410):398--409, 1990.

\bibitem{ghahramani2001propagation}
Zoubin Ghahramani and Matthew~J Beal.
\newblock Propagation algorithms for variational bayesian learning.
\newblock In {\em Advances in neural information processing systems}, pages
  507--513, 2001.

\bibitem{ghosal2000}
Subhashis Ghosal, Jayanta~K. Ghosh, and Aad~W. van~der Vaart.
\newblock Convergence rates of posterior distributions.
\newblock {\em Ann. Statist.}, 28(2):500--531, 2000.

\bibitem{ghosal2001entropies}
Subhashis Ghosal and Aad~W van~der Vaart.
\newblock Entropies and rates of convergence for maximum likelihood and bayes
  estimation for mixtures of normal densities.
\newblock {\em Annals of Statistics}, pages 1233--1263, 2001.

\bibitem{graves2011practical}
Alex Graves.
\newblock Practical variational inference for neural networks.
\newblock In {\em Advances in Neural Information Processing Systems}, pages
  2348--2356, 2011.

\bibitem{ho2016strong}
Nhat Ho, XuanLong Nguyen, et~al.
\newblock On strong identifiability and convergence rates of parameter
  estimation in finite mixtures.
\newblock {\em Electronic Journal of Statistics}, 10(1):271--307, 2016.

\bibitem{hoffman2013stochastic}
Matthew~D Hoffman, David~M Blei, Chong Wang, and John Paisley.
\newblock Stochastic variational inference.
\newblock {\em The Journal of Machine Learning Research}, 14(1):1303--1347,
  2013.

\bibitem{humphreys2000approximate}
K~Humphreys and DM~Titterington.
\newblock {Approximate Bayesian inference for simple mixtures}.
\newblock In {\em Proc. Computational Statistics 2000}, pages 331--336, 2000.

\bibitem{jordan1999introduction}
Michael~I Jordan, Zoubin Ghahramani, Tommi~S Jaakkola, and Lawrence~K Saul.
\newblock An introduction to variational methods for graphical models.
\newblock {\em Machine learning}, 37(2):183--233, 1999.

\bibitem{kingma2014adam}
Diederik Kingma and Jimmy Ba.
\newblock Adam: A method for stochastic optimization.
\newblock {\em arXiv preprint arXiv:1412.6980}, 2014.

\bibitem{kushner1997stochastic}
Harold~J Kushner and G~George Yin.
\newblock Stochastic approximation algorithms and applications.
\newblock 1997.

\bibitem{lecam1973convergence}
Lucien LeCam.
\newblock Convergence of estimates under dimensionality restrictions.
\newblock {\em The Annals of Statistics}, pages 38--53, 1973.

\bibitem{mackay1997ensemble}
David~JC MacKay.
\newblock Ensemble learning for hidden markov models.

\bibitem{nguyen2013convergence}
XuanLong Nguyen.
\newblock Convergence of latent mixing measures in finite and infinite mixture
  models.
\newblock {\em The Annals of Statistics}, 41(1):370--400, 2013.

\bibitem{pati2015optimal}
Debdeep Pati, Anirban Bhattacharya, and Guang Cheng.
\newblock {Optimal Bayesian estimation in random covariate design with a
  rescaled Gaussian process prior}.
\newblock {\em Journal of Machine Learning Research}, 16:2837--2851, 2015.

\bibitem{pati2014posterior}
Debdeep Pati, Anirban Bhattacharya, Natesh~S Pillai, and David Dunson.
\newblock Posterior contraction in sparse bayesian factor models for massive
  covariance matrices.
\newblock {\em The Annals of Statistics}, 42(3):1102--1130, 2014.

\bibitem{pati2013posterior}
Debdeep Pati, David~B Dunson, and Surya~T Tokdar.
\newblock Posterior consistency in conditional distribution estimation.
\newblock {\em Journal of multivariate analysis}, 116:456--472, 2013.

\bibitem{rousseau2016frequentist}
Judith Rousseau.
\newblock {On the frequentist properties of Bayesian nonparametric methods}.
\newblock {\em Annual Review of Statistics and Its Application}, 3:211--231,
  2016.

\bibitem{rousseau2017asymptotic}
Judith Rousseau and Botond Szabo.
\newblock Asymptotic behaviour of the empirical bayes posteriors associated to
  maximum marginal likelihood estimator.
\newblock {\em The Annals of Statistics}, 45(2):833--865, 2017.

\bibitem{schwartz1965bayes}
Lorraine Schwartz.
\newblock {On Bayes procedures}.
\newblock {\em Probability Theory and Related Fields}, 4(1):10--26, 1965.

\bibitem{stephens2000dealing}
Matthew Stephens.
\newblock Dealing with label switching in mixture models.
\newblock {\em Journal of the Royal Statistical Society: Series B (Statistical
  Methodology)}, 62(4):795--809, 2000.

\bibitem{ueda2002bayesian}
Naonori Ueda and Zoubin Ghahramani.
\newblock Bayesian model search for mixture models based on optimizing
  variational bounds.
\newblock {\em Neural Networks}, 15(10):1223--1241, 2002.

\bibitem{van2000asymptotic}
Aad~W Van~der Vaart.
\newblock {\em Asymptotic statistics}, volume~3.
\newblock Cambridge university press, 2000.

\bibitem{wainwright2008graphical}
Martin~J Wainwright and Michael~I Jordan.
\newblock Graphical models, exponential families, and variational inference.
\newblock {\em Foundations and Trends{\textregistered} in Machine Learning},
  1(1--2):1--305, 2008.

\bibitem{wang2005inadequacy}
Bo~Wang and DM~Titterington.
\newblock {Inadequacy of interval estimates corresponding to variational
  Bayesian approximations}.
\newblock In {\em AISTATS}, 2005.

\bibitem{yang2014minimax}
Yun Yang and David~B Dunson.
\newblock {Minimax optimal Bayesian aggregation}.
\newblock {\em arXiv preprint arXiv:1403.1345}, 2014.

\bibitem{yang2017variational}
Yun Yang, Debdeep Pati, and Anirban Bhattacharya.
\newblock $\alpha$-variational inference with statistical guarantees.
\newblock {\em arXiv preprint arXiv:1710.03266}, 2017.

\end{thebibliography}

\end{document}